\documentclass[automark,paper=a4,fontsize=10pt]{scrartcl}

\usepackage{ae}
\usepackage{amsmath}
\usepackage{amssymb}
\usepackage{amsthm}
\usepackage{amsfonts}
\usepackage{dsfont}
\usepackage[mathscr]{eucal}
\usepackage{enumerate}
\usepackage[normalem]{ulem}
\usepackage{scrpage2}
\usepackage{pstricks,pst-node,pst-tree}
\usepackage{arydshln}
\usepackage{rotating}

\setheadsepline{0.4pt}

\ihead{On associative spectra of operations}
\ohead{\headmark}
\newtheorem{thm}{Theorem}[section]
\newtheorem{pro}[thm]{Proposition}
\newtheorem{defn}[thm]{Definition}
\newtheorem{cor}[thm]{Corollary}
\newtheorem{ex}[thm]{Example}
\newtheorem{rem}[thm]{Remark}
\newtheorem{lemma}[thm]{Lemma}

\newcommand{\bthm}{\begin{thm}\rm}
\newcommand{\bpro}{\begin{pro}\rm}
\newcommand{\bd}{\begin{defn}\rm}
\newcommand{\brem}{\begin{rem}\rm}
\newcommand{\bcor}{\begin{cor}\rm}
\newcommand{\bex}{\begin{ex}\rm}
\newcommand{\bl}{\begin{lemma}\rm}
\newcommand{\ethm}{\end{thm}\noindent}
\newcommand{\epro}{\end{pro}\noindent}
\newcommand{\ed}{\end{defn}\noindent}
\newcommand{\erem}{\end{rem}\noindent}
\newcommand{\ecor}{\end{cor}\noindent}
\newcommand{\eex}{\end{ex}\noindent}
\newcommand{\el}{\end{lemma}\noindent}

\newcommand{\bproof}{\textbf{Proof: }}
\newcommand{\eproof}{\nopagebreak \vspace{-2.0ex} \begin{flushright} \tiny $\blacksquare$ 
                                                \end{flushright} \smallskip}

\newcommand{\N}{\mathds{N}}
\newcommand{\Nplus}{\mathds{N^+}}

\newcommand{\Z}{\mathds{Z}}

\newcommand{\C}{\mathds{C}}

\newcommand{\ds}{\displaystyle}
\renewcommand{\Tp}{{\bf T^{(p)}}}
\newcommand{\T}{T_{\omega} (x)}
\newcommand{\TX}{T_{\omega} (X)}
\newcommand{\w}{\omega^{{\bf T^{(p)}}}}
\newcommand{\occN}[1]{\left| #1 \right|_{\omega}}
\newcommand{\length}[1]{\left| #1 \right|}
\renewcommand{\l}[1]{\ell^{(#1)}}

\newcommand{\B}[2]{B_{#1}^{(#2)}}
\newcommand{\lset}[2]{\left\{ \left. #1 \ \right| \ #2 \right\}}
\newcommand{\rset}[2]{\left\{ #1 \ \left| \ #2 \right. \right\}}
\newcommand{\set}[1]{\left\{ #1 \right\}}

\newcommand{\en}[1]{\varepsilon_{#1}}
\newcommand{\function}[5]{\begin{array}{cccc}
   #1: & #2 & \longrightarrow & #3 \\
       & #4 & \longmapsto     & #5
\end{array}}
\newcommand{\functionhead}[3]{#1: #2 \longrightarrow #3}

\renewcommand{\C}[2]{C_{#1}^{(#2)}}
\renewcommand{\L}[1]{\Lambda^{(#1)}}
\newcommand{\pair}[2]{\left( #1, #2 \right)}

\newcommand{\Alg}[1]{\textrm{Alg}\left( #1 \right)}

\DeclareMathOperator{\Mod}{Mod}
\DeclareMathOperator{\Id}{Id}
\DeclareMathOperator{\ModBrack}{Mod_{Brack}}
\DeclareMathOperator{\IdBrack}{Id_{Brack}}
\newcommand{\modelsBrack}{\models_{\textrm{Brack}}}
\newcommand{\eq}[2]{#1 \approx #2}
\newcommand{\K}{{\mathscr K}}

\newcommand{\fS}[1]{\sigma \left( #1 \right)}
\newcommand{\fSn}[2]{\sigma_{#2} \left( #1 \right)}
\newcommand{\aS}[2]{s_{#2} \left( #1 \right)}

\DeclareMathOperator{\Eq}{Eq}
\newcommand{\TpX}{{\bf T^{(p)}(X)}}
\newcommand{\Ec}[1]{\left[ #1 \right]}
\DeclareMathOperator{\ST}{IT}
\newcommand{\M}[3]{M_{#1,#2}^{(#3)}}

\newcommand{\p}[1]{{\mathcal P} \left( #1 \right)}
\newcommand{\FSB}{{\bf FS}}
\newcommand{\FS}{FS}

\newcommand{\myfootnote}[1]{
\renewcommand{\thefootnote}{}
\footnotetext{#1}
\renewcommand{\thefootnote}{\arabic{footnote}}}

\begin{document}

\date{}
\author{Sebastian Liebscher\footnote{TU Dresden, Institut f\"ur Algebra, Dresden, Germany; e-mail: Seb.Liebscher@gmx.de}~~and Tam\'as
Waldhauser\footnote{University of Szeged, Bolyai Institute, Szeged, Hungary; e-mail: twaldha@math.u-szeged.hu}}
\title{On associative spectra of operations}
\maketitle

\noindent
{\small
{\bf Abstract.}
The distance of an operation from being associative can be ``measured'' by its
associative spectrum, an appropriate sequence of positive integers. Associative
spectra were introduced in a publication by B. Cs\'ak\'any and T. Waldhauser in
2000 for binary operations (see \cite{CsakanyWaldhauser}).
We generalize this concept to $2 \le p$-ary operations,
interpret associative spectra in terms of equational theories, and use this interpretation
to find a characterization of fine spectra, to construct polynomial associative spectra, and to
show that there are continuously many different spectra. Furthermore, an equivalent representation
of bracketings is studied.}

\section{Introduction}

\myfootnote{Research supported by the Hungarian National Foundation for Scientific Research grant no. T48809 and K60148.}
B. Cs\'ak\'any and T. Waldhauser introduced associative spectra for binary operations in\cite{CsakanyWaldhauser}.
The main focus point in their paper was the spectrum of groupoids with two or three elements.

In this paper, we generalize first in Section 2 the definition of associative spectrum to $2 \le p$-ary operations
with the help of special unary terms, which will be called bracketings. Enumerations are used to distinguish between the
variable symbols in a bracketing. Using these enumarations it is possible to define a reduct $\ModBrack - \IdBrack$ of the well-known
{\sc Galois}-connection $\Mod - \Id$. The {\sc Galois}-closed sets on the side of the identities are called fine spectra, which is a
refinement of the notion of associative spectrum. Finally, some useful operations on bracketings are defined which are needed in the characterization
of fine spectra.
In Section 3 we give the characterization of fine spectra and a first application of it, a generalization of the generalized associative law.
After that, insertion tuples are developed as an equivalent representation of bracketings in Section 4. With the help of these tuples
the explicit formula of the generalized {\sc Catalan} numbers is proven, where the generalized {\sc Catalan} numbers
count bracketings of a given length.
In Section 5 three different polynomial spectra are presented which solve problem 3 in \cite{CsakanyWaldhauser}.
The lattice of fine spectra is studied in Section 6. The covering relation, the atoms and the coatoms of this lattice are described.
Furthermore, it is shown that there are continuously many different spectra.
In Section 7 we look at some examples of finite groupoids (one of them has a polynomial spectrum from Section 5). It is shown that
every finally associative spectrum appears as the fine spectrum of a finite groupoid.
We conclude in Section 8 with the formulation of a few open problems.

\section{Definitions and notation}\label{sectionOne}

The {\it algebra of $p$-ary bracketings} is defined as the term algebra
\[
 \Tp := \pair{\T}{\w},
\]
where $p$ is a natural number greater or equal to 2, the alphabet is $\{ x \}$ and the signature is $\{ \omega \}$ with $\omega$ as a
$p$-ary operation symbol.

We call the (unary) terms $t \in \T$ {\it $p$-ary bracketings} or simply bracketings if $p$ is known.
The {\it occurence number} $\occN{t}$ of a bracketing $t \in \T$ is defined as the number of occurences of the
operation symbol $\omega$ in $t$. The following trivial equalities hold:
\[
 \occN{x} = 0 \ \textrm{ and } \ \forall t_1, \ldots, t_p \in \T : \ \occN{\omega t_1 \ldots t_p} = 1 + \sum_{k=1}^{p} \occN{t_k}.
\]
The {\it length} $\length{t}$ of a bracketing $t \in \T$ is defined as the number of occurences of the
variable symbol $x$ in $t$. It is an easy observation to show that $\length{t} = (p - 1) \cdot \occN{t} + 1$ holds for all
bracketings $t \in \T$. The length could be defined recursively too:
\[
 \length{x} = 1 \ \textrm{ and } \ \forall t_1, \ldots, t_p \in \T : \ \length{\omega t_1 \ldots t_p} = \sum_{k=1}^{p} \length{t_k}.
\]
Bracketings with occurence number $n$ can be viewed as trees with branching factor $p$, $n$ inner nodes and $(p-1) \cdot n + 1$ leafs
(because the symbols $\omega$ are the inner nodes and the symbols $x$ are the leafs).
We denote by
\[
 \B{n}{p} := \lset{t \in \T}{\occN{t} = n}
\]
the set of all bracketings with occurence number $n$. The {\it length function}, which transforms occurence numbers into lengths, is given by
\[
 \function{\l{p}}{\N}{\Nplus}{n}{(p - 1) \cdot n + 1.}
\]
For example, the first binary bracketings are given in Table \ref{TableOne} (insertion tuples are defined in Definition~\ref{DefInsertionTuple}).

\begin{sidewaystable}
\centering
\renewcommand{\arraystretch}{1.5}
\newcommand{\bonusDepth}{5pt}
\newcommand{\bonusWidth}{19.6pt}
\setlength{\unitlength}{1mm}
\psset{nodesep=2pt,levelsep=20pt,treesep=15pt}

\noindent
\begin{tabular}{|l|c|c|c|c|}
 \hline
 & $\B{0}{2}$ & $\B{1}{2}$ & \multicolumn{2}{|c|}{$\B{2}{2}$} \\ \hline
 tree &
 \TR{$x$} &
 \pstree[xbbd=\bonusDepth]{\TR{$\omega$}} {
        \TR{$x$}
        \TR{$x$}
 } &
\pstree[xbbd=\bonusDepth]{\TR{$\omega$}} {
    \pstree{\TR{$\omega$}}
    {
        \TR{$x$}
        \TR{$x$}
    }
    \TR{$x$}
} &
\pstree[xbbd=\bonusDepth]{\TR{$\omega$}} {
    \TR{$x$}
    \pstree{\TR{$\omega$}}
    {
        \TR{$x$}
        \TR{$x$}
    }
} \\ \hline
 bracketing & \hspace{\bonusWidth} $x$ \hspace{\bonusWidth}
            & \hspace{\bonusWidth} $\omega x x$ \hspace{\bonusWidth}
            & \hspace{\bonusWidth} $\omega \omega xxx$ \hspace{\bonusWidth}
            & \hspace{\bonusWidth} $\omega x \omega xx$ \hspace{\bonusWidth} \\ \hline
 bracketing (infix) & $x$ & $(xx)$ & $((xx)x)$ & $(x(xx))$ \\ \hline
 insertion tuple & $\emptyset$ & $(1)$ & $(1,1)$ & $(1,2)$ \\ \hline
\end{tabular}\bigskip

\noindent
\begin{tabular}{|l|c|c|c|c|c|}
 \hline
 & \multicolumn{5}{|c|}{$\B{3}{2}$} \\ \hline
 tree & \pstree[xbbd=\bonusDepth]{\TR{$\omega$}} {
    \pstree{\TR{$\omega$}}
    {
        \pstree{\TR{$\omega$}}
        {
            \TR{$x$}
            \TR{$x$}
        }
        \TR{$x$}
    }
    \TR{$x$}
} &
\pstree[xbbd=\bonusDepth]{\TR{$\omega$}} {
    \pstree{\TR{$\omega$}}
    {
        \TR{$x$}
        \pstree{\TR{$\omega$}}
        {
            \TR{$x$}
            \TR{$x$}
        }
    }
    \TR{$x$}
} &
\pstree[xbbd=\bonusDepth]{\TR{$\omega$}} {
    \pstree{\TR{$\omega$}}
    {
        \TR{$x$}
        \TR{$x$}
    }
    \pstree{\TR{$\omega$}}
    {
        \TR{$x$}
        \TR{$x$}
    }
} &
\pstree[xbbd=\bonusDepth]{\TR{$\omega$}} {
    \TR{$x$}
    \pstree{\TR{$\omega$}}
    {
        \pstree{\TR{$\omega$}}
        {
            \TR{$x$}
            \TR{$x$}
        }
        \TR{$x$}
    }
} &
\pstree[xbbd=\bonusDepth]{\TR{$\omega$}} {
    \TR{$x$}
    \pstree{\TR{$\omega$}}
    {
        \TR{$x$}
        \pstree{\TR{$\omega$}}
        {
            \TR{$x$}
            \TR{$x$}
        }
    }
} \\ \hline
 bracketing & $\omega \omega \omega x x x x$ & $\omega \omega x \omega x x x$ & $\omega \omega x x \omega x x$ & $\omega x \omega \omega x x x$
 & $\omega x \omega x \omega x x$ \\ \hline
 bracketing (infix) & $(((xx)x)x)$ & $((x(xx))x)$ & $((xx)(xx))$ & $(x((xx)x))$ & $(x(x(xx)))$ \\ \hline
 insertion tuple & $(1,1,1)$ & $(1,1,2)$ & $(1,1,3)$ & $(1,2,2)$ & $(1,2,3)$ \\ \hline
\end{tabular}
\caption{Binary bracketings, their tree correspondences and their insertion tuples}\label{TableOne}
\end{sidewaystable}
In the next step we want to distinguish between the variable symbols $x$ in a bracketing. Therefore, we define the
{\it enumerations $\functionhead{\en{j}}{\T}{\TX}$} by term induction as follows,
 where $\TpX = \pair{\TX}{\omega^{\TpX}}$ is the term algebra over the alphabet $X = \rset{x_i}{i \in \Nplus}$:
\begin{itemize}
 \item $\forall j \in \Nplus : \ \en{j} \left( x \right) = x_j;$
 \item $\forall t_1, \ldots, t_p \in \T \ \forall j \in \Nplus : \ \en{j} \left( \omega t_1 \ldots t_p \right) = \omega
        \en{j_1} \left( t_1 \right) \ldots \en{j_p} \left( t_p \right)$ with $\displaystyle j_m := j + \sum_{k=1}^{m - 1} \length{t_k}$.
\end{itemize}
It is obvious that $\en{j} \left( t \right)$ contains exactly the variable symbols $\set{x_j, \ldots, x_{j + \length{t} - 1}}$.
As an example we look again at some binary bracketings:
\[
 \en{j} (\omega \omega x x x ) = \omega \omega x_j x_{j+1} x_{j+2} \ \textrm{, } \
 \en{j} (\omega x \omega x x ) = \omega x_j \omega x_{j+1} x_{j+2}.
\]
For a simpler notation we denote by
\[
 \L{p} := \lset{\pair{s}{t} \in \T \times \T}{\occN{s} = \occN{t}}
\]
the relation of all bracketings with the same occurence number.
It is an easy observation that $\L{p}$ is a congruence relation of $\Tp$.
Further on, we will denote pairs $\pair{s}{t} \in \L{p}$ simply by $\eq{s}{t}$ and we will call them {\it identities}.
From the example above we know that these identities can be interpreted via enumaration as generalized associativity conditions. \bigskip

\noindent
We call an algebra ${\bf A}$ to the signature $\{ \omega \}$ a $p$-ary groupoid and denote it by
${\bf A} \in \Alg{\omega}$.
Now we can define a reduct of the well-known {\sc Galois}-connection $\Mod-\Id$.
Let $\modelsBrack \subseteq \Alg{\omega} \times \L{p}$ be defined as
\[
 {\bf A} \modelsBrack \eq{s}{t} :\Longleftrightarrow {\bf A} \models \pair{\en{1} \left( s \right)}{\en{1} \left( t \right)}.
\]
Because of the full invariance of $\Id {\bf A}$ it is obvious that
\[
 {\bf A} \modelsBrack \eq{s}{t} \Longrightarrow \forall j \in \Nplus : \ {\bf A} \models \pair{\en{j} \left( s \right)}{\en{j} \left( t \right)}.
\]
The {\sc Galois}-closed sets are given for any $\Sigma \subseteq \L{p}$ and $\K \subseteq \Alg{\omega}$ by
\begin{itemize}
 \item $\ModBrack \Sigma := \rset{{\bf A} \in \Alg{\omega}}{\forall \eq{s}{t} \in \Sigma : \ {\bf A} \modelsBrack \eq{s}{t}}$,

       which is of course a special variety that has additional properties (see open problems);
 \item $\IdBrack \K := \lset{\eq{s}{t} \in \L{p}}{\forall {\bf A} \in \K : \ {\bf A} \modelsBrack \eq{s}{t}}$,

       which is a reduct of the equational theory of $\K$.
\end{itemize}
We will further on denote any $\Sigma \subseteq \L{p}$ equivalently by the sequence
\[
 \left( \Sigma_n \right)_{n \in \N} := \left( \Sigma \cap \left( \B{n}{p} \times \B{n}{p} \right) \right)_{n \in \N}.
\]
For a $p$-ary groupoid ${\bf A}$ we define two different spectra:
\begin{itemize}
 \item the \emph{fine spectrum} of ${\bf A}$: $\fS{{\bf A}} := \IdBrack{{\bf A}}$, or equivalently (see above)
       $\left( \fSn{{\bf A}}{n} \right)_{n \in \N}$ with $\fSn{{\bf A}}{n} = \fS{{\bf A}} \cap \left( \B{n}{p} \times \B{n}{p} \right)$;
 \item the \emph{associative spectrum} of ${\bf A}$: $\left( \aS{{\bf A}}{n} \right)_{n \in \N}
       := \left( \left| \B{n}{p} / \fSn{{\bf A}}{n} \right| \right)_{n \in \N}$.
\end{itemize}
We say that ${\bf A}$ is {\it associative} iff $\aS{{\bf A}}{2} = 1$. The following two observations are trivial.
\bpro
If ${\bf A} \in \Alg{\omega}$ is a subgroupoid or a homomorphic image of ${\bf B} \in \Alg{\omega}$, then
\[
 \fS{{\bf A}} \supseteq \fS{{\bf B}} \ \textrm{ and } \ \forall n \in \N : \ \aS{{\bf A}}{n} \le \aS{{\bf B}}{n}.
\]
\epro

\bpro
If ${\bf A} \in \Alg{\omega}$ and ${\bf B} \in \Alg{\omega}$ are isomorphic or antiisomorphic, then their spectra coincide:
\[
 \fS{{\bf A}} = \fS{{\bf B}}.
\]
\epro
Finally, we define some useful operations for bracketings:
\begin{itemize}
 \item $\functionhead{\gamma_i}{\T}{\T}$ ($i = 1, \ldots, p$) is defined as
       \[
        \function{\gamma_i}{\T}{\T}{t}{\omega t_1 \ldots t_p}
       \]
       with $t_i = t$ and $t_k = x$ for all $k \in \set{1, \ldots, p} \setminus \set{i}$. So $\gamma_i (t) = \omega x \ldots x \, t \
       x \ldots x$ is the insertion of
       $t$ at the $i$-th position in $\omega x \ldots x$.
 \item For the definition of $\beta_i$ ($i \in \Nplus$) we need some auxiliary functions $\alpha_i$:
       \[
        \function{\alpha_i}{X}{\T}{x_k}{\begin{cases}
                                         \omega x \ldots x \in \B{1}{p}, & \textrm{if } k = i; \\
                                         x,                              & \textrm{otherwise}.
                                        \end{cases}}
       \]
       Denote (here and further on)
       by $\functionhead{\alpha_i^{\#}}{\TX}{\T}$ the unique homomorphism that continues $\alpha_i$. Then $\beta_i$ is defined as
       \[
        \function{\beta_i}{\T}{\T}{t}{\alpha_i^{\#} \left( \en{1} \left( t \right) \right).}
       \]
       So $\beta_i (t)$ is the insertion of $\omega x \ldots x$ at the $i$-th symbol $x$  in $t$ (if present).
\end{itemize}
It is easy to check that for any bracketing $t \in \B{n}{p}$ the resulting bracketings $\gamma_i \left( t \right)$ ($i = 1, \ldots, p$)
and $\beta_i \left( t \right)$ ($i = 1, \ldots, \l{p} (n)$) are in $\B{n+1}{p}$. \bigskip

\noindent
To put these operators together we define for any positive natural number $n \in \N$ the {\it implication operator}
$\delta_n$ as follows:
\[
 \function{\delta_n}{\Eq \B{n}{p}}{\Eq \B{n+1}{p}}{\pi}{\displaystyle \left( \bigcup_{\xi \in \set{\gamma_1, \ldots, \gamma_p,
 \beta_1, \ldots, \beta_{\l{p} \left( n \right)}}} \lset{\eq{\xi \left( s \right)}{\xi \left( t \right)}}{\eq{s}{t} \in \pi} \right)^{*},}
\]
where $\Eq \B{n}{p}$ denotes the set of equivalence relations on $\B{n}{p}$ and $\tau^{*}$ denotes the transitive closure of
$\tau$.

\section{Characterization of fine spectra}

Our main goal is to characterize the {\sc Galois}-closed sets $\IdBrack \K$ and the fine spectra of arbitrary groupoids. \bigskip

\noindent
First we need three preparatory lemmata. The first one shows a recursion formula for the operators $\beta_i$.
\bl\label{LemmaBeta}
For all $k \in \set{1, \ldots, p}$ and for all $t_1, \ldots, t_p \in \T$ with $i \in \set{1, \ldots, \length{t_k}}$ we have:
\[
 \omega \, t_1 \ldots t_{k-1} \, \beta_i \left( t_k \right) \, t_{k+1} \ldots t_p = \beta_j \left( \omega t_1 \ldots t_p \right),
\]
where $\displaystyle j := i + \sum_{l=1}^{k-1} \length{t_l}$.
\el
The proof is left to the reader; it is just a transformation of the insertion index.
%
The next statement shows that all bracketings can be obtained with the operators $\beta_i$ starting with $x$.
\bl\label{LemmaBetaConstructsEverything}
For all $n \in \N$ we have:
\[
 \B{n+1}{p} = \rset{\beta_i \left( t \right)}{t \in \B{n}{p}, i = 1, \ldots, \l{p} (n)}.
\]
\el
This follows directly with the previous lemma by induction on $n$.

\bl\label{LemmaClosedUnderSubstitutions}
If $\Sigma \subseteq \L{p}$ is an equivalence relation that is closed under the implication operator, i.e.
\[
 \forall n \in \N : \ \delta_n \left( \Sigma_n \right) \subseteq \Sigma_{n+1},
\]
then
\[
 \eq{s}{t} \in \Sigma \Longrightarrow \eq{\alpha^{\#} \left( \en{1} (s) \right)}{\alpha^{\#} \left( \en{1} (t) \right)} \in \Sigma
\]
holds for all $\functionhead{\alpha}{X}{\T}$.
\el
\bproof
We choose an arbitrary but fixed identity $\eq{s}{t} \in \Sigma$.
Then we apply induction on $\displaystyle n:=\sum_{i=1}^{\length{s}} \occN{\alpha \left( x_i \right)}$:
 For $n=0$, $\functionhead{\alpha}{X}{\T}$ must
               map each $x_i$ on $x$ ($i = 1, \ldots, \length{s}$), thus
               \[
                \eq{\alpha^{\#} \left( \en{1} (s) \right)}{\alpha^{\#} \left( \en{1} (t) \right)} =
                \eq{s}{t} \in \Sigma.
               \]
 For the induction step from $n$ to $n+1$ let $k \in \set{1, \ldots, \length{s}}$ be a position where the occurence number is greater than 0, i.e.
                             $\occN{\alpha \left( x_k \right)} > 0$. With the previous Lemma~\ref{LemmaBetaConstructsEverything} we can
                             find a bracketing $t_k$ and a natural number
                             $j \in \set{1, \ldots, \length{t_k}}$
                             such that
                             \[
                              \beta_j \left( t_k \right) = \alpha \left( x_k \right).
                             \]
                             Now we can define a reduct of $\alpha$:
                             \[
                              \function{\tilde{\alpha}}{X}{\T}{x_i}{\begin{cases}
                                                                   \alpha \left( x_i \right), & \textrm{if } i \ne k; \\
                                                                   t_k,                       & \textrm{if } i = k.
                                                                  \end{cases}}
                             \]
                             It is easy to see that $\displaystyle \sum_{i=1}^{\length{s}} \occN{\tilde{\alpha} \left( x_i \right)} = n$ holds.
                             With $\displaystyle l := j + \sum_{m=1}^{k-1} \length{\alpha \left( x_m \right)}$
                             it can be shown with the help of Lemma \ref{LemmaBeta} that
                             \[
                              \beta_{l} \circ \tilde{\alpha}^{\#} \circ \en{1} = \alpha^{\#} \circ \en{1}
                             \]
                             holds.
                             Then together with the induction hypothesis for $\tilde{\alpha}^{\#}$
                             and the prerequisite that $\Sigma$
                             is closed under the implication operator we get
               \[
                \eq{\alpha^{\#} \left( \en{1} (s) \right)}{\alpha^{\#} \left( \en{1} (t) \right)} \in \Sigma.
               \]
\eproof
Now we are able to prove our main result, which is an analogon of the well-known characterization of equational theories.
\bthm\label{TheoremFS}
For any $\K \subseteq \Alg{\omega}$ and $\Sigma \subseteq \L{p}$ the following hold:
\begin{enumerate}[(a)]
 \item If $\Sigma$ is an equivalence relation that is closed under the implication operator
       then $\Sigma$ is a congruence of $\Tp$.
 \item $\IdBrack \K$ is closed under the implication operator.
 \item If $\Sigma$ is an equivalence relation that is closed under the implication operator then
       \[
        \IdBrack \set{\Tp / \Sigma} = \Sigma.
       \]
\end{enumerate}
\ethm
\bproof
\begin{enumerate}[(a)]
 \item It suffices to show that for all $i \in \set{1, \ldots, p}$ , $\eq{s}{t} \in \Sigma$ and $t_1, \ldots, t_p \in \T$ it holds:
       \[
        \eq{\omega t_1 \ldots t_{i-1} \, s \, t_{i+1} \ldots t_p}{\omega t_1 \ldots t_{i-1} \, t \, t_{i+1} \ldots t_p} \in \Sigma.
       \]
       (The general case follows then by applying this rule repeatedly on each position). We know that
       $\Sigma$ is closed under the implication operator, so:
       \[
        \eq{\gamma_i (s)}{\gamma_i (t)} = \eq{\omega x \ldots x \, s \, x \ldots x}{\omega x \ldots x \, t \, x \ldots x} \in \Sigma.
       \]
       According to Lemma~\ref{LemmaBetaConstructsEverything} we have a sequence of natural numbers
       $i_1, \ldots, i_{\occN{t_p}} $ such that
       \[
        t_p = \beta_{i_{\occN{t_p}}} \circ \ldots \circ \beta_{i_1} \left( x \right).
       \]
       So it follows from Lemma~\ref{LemmaBeta} with $j_k := i_k + p - 2 + \length{s}$ that
       \[
        \begin{array}{rcl}
         \beta_{j_{\occN{t_p}}} \circ \ldots \circ \beta_{j_1} \left( \gamma_i (s) \right) & = &
         \omega x \ldots x \, s \, x \ldots x \, \left( \beta_{i_{\occN{t_p}}} \circ \ldots \circ \beta_{i_1} \left( x \right) \right) \\
         & = & \omega x \ldots x \, s \, x \ldots x \, t_p
        \end{array}
       \]
       holds. Similarly we have $\beta_{j_{\occN{t_p}}} \circ \ldots \circ \beta_{j_1} \left( \gamma_i (t) \right) =
        \omega x \ldots x \, t \, x \ldots x \, t_p$ and
       since $\Sigma$ is closed under the implication operator, we get
       \[
        \eq{\omega x \ldots x \, s \, x \ldots x \, t_p}{\omega x \ldots x \, t \, x \ldots x \, t_p} \in \Sigma.
       \]
       This construction step can be repeated, and finally
       we obtain
       \[
        \eq{\omega t_1 \ldots t_{i-1} \, s \, t_{i+1} \ldots t_p}{\omega t_1 \ldots t_{i-1} \, t \, t_{i+1} \ldots t_p} \in \Sigma.
       \]
 \item This result follows obviously from the fact that $\Id \K$ is always a fully invariant congruence relation of $\TX$
       with the help of the following two equations:
       \[
         \begin{array}{rcl}
          \ds \en{1} \left( \gamma_i \left( t \right) \right) & \ds = & \ds \omega x_1 \ldots x_{i - 1} \, \en{i} (t) \, x_{i+\length{t}}
          \ldots x_{p-1+\length{t}}, \smallskip \\
          \ds \en{1} \left( \beta_i \left( t \right) \right) & \ds = & \ds \tilde{\alpha}_i^{\#} \left( \en{1} (t) \right),
         \end{array}
       \]
       where
       \[
        \function{\tilde{\alpha}_i}{X}{\TX}{x_k}{\begin{cases}
                                         x_k,                         & \textrm{if } k < i; \\
                                         \omega x_k \ldots x_{k+p-1}, & \textrm{if } k = i; \\
                                         x_{k+p-1},                   & \textrm{if } k > i.
                                        \end{cases}}
       \]
 \item
   \begin{enumerate}
    \item[``$\subseteq$'':] Let $\eq{s}{t} \in \IdBrack \set{\Tp / \Sigma}$. Then we have $\Tp / \Sigma \models \eq{\en{1} (s)}{\en{1} (t)}$.
                            With the full invariance of $\Id \set{\Tp / \Sigma}$ and the function
                            \[
                             \function{\alpha}{X}{\Tp / \Sigma}{x_k}{\Ec{x}}
                            \]
                            that maps each $x_k$ to the equivalence class of $x$ in $\T / \Sigma$
                            it follows that $\Ec{s} = \alpha^{\#} \left( \en{1} (s) \right) = \alpha^{\#} \left( \en{1} (t) \right)
                             = \Ec{t}$.
                            Therefore, $\eq{s}{t}$ is in $\Sigma$.
    \item[``$\supseteq$'':] Follows directly with Lemma~\ref{LemmaClosedUnderSubstitutions}.
   \end{enumerate}
\end{enumerate}
\eproof
The following corollary summarizes the last theorem.
\bcor\label{TheoremCharFS}
For $\Sigma \subseteq \L{p}$ the following are equivalent:
\begin{enumerate}[(a)]
 \item $\Sigma$ is an equivalence relation and $\Sigma$ is closed under the implication operator.
 \item $\IdBrack \ModBrack \Sigma = \Sigma$.
 \item There exists a groupoid {\bf A} such that $\fS{{\bf A}} = \Sigma$.
\end{enumerate}
\ecor
As a first application of our main result we show that the general associative law and a generalization of it hold.
\bthm\label{TheoremGenAssLaw}
For any $p$-ary groupoid ${\bf A}$ the following hold:
\begin{enumerate}[(a)]
 \item $\aS{{\bf A}}{0} = \aS{{\bf A}}{1} = 1$;
 \item\label{GenAssLawB} ${\bf A}$ is associative $\Longleftrightarrow \forall n \in \N : \ \aS{{\bf A}}{n} = 1$;
 \item $\aS{{\bf A}}{n} = 1 \Longrightarrow \forall m \in \N, m \ge n : \ \aS{{\bf A}}{m} = 1$ for any $n \ge 2$.
\end{enumerate}
If the conclusion of (c) is fulfilled then we say that ${\bf A}$ and the associative spectrum of ${\bf A}$ are {\it finally associative}.
\ethm
\bproof
\begin{enumerate}[(a)]
 \item Absolutely clear because $\B{0}{p} = \set{x}$ and $\B{1}{p} = \set{\omega x \ldots x}$.
 \item The direction ``$\Longleftarrow$'' is clear from the definition in Section~\ref{sectionOne}.
       For the other direction we suppose that ${\bf A}$ is associative. For an arbitrary natural number $n \ge 2$ let
       $t \in \B{n}{p}$
       be the left associated bracketing $\omega \ldots \omega x \ldots x$, which means that all symbols $\omega$ occur before the $x$'s.
       Assume we have another bracketing $s \in \B{n}{p}$. We will prove that ${\bf A} \modelsBrack \eq{s}{t}$.
       Since $s$ is not left associated, it has a subbracketing of the form
       \[
         \omega \left( t_1 \ldots t_k \, \omega \left( s_1 \ldots s_p \right) \ldots \right),
       \]
       where $k \ge 1$. Then we have by associativity
       \[
        {\bf A} \modelsBrack \eq{\omega \left( t_1 \ldots t_k \, \omega \left( s_1 \ldots s_p \right) \ldots \right)} 
        {\omega \left( \omega \left( t_1 \ldots t_k \, s_1 \ldots s_{p-k} \right) s_{p-k+1} \ldots s_p \ldots \right)}.
       \]
       This way one $\omega$ is moving to the left, and after finitely
       many such steps we reach $t$, i.e. ${\bf A} \modelsBrack \eq{s}{t}$.
 \item We show this fact via contradiction: Assume that there exists a groupoid ${\bf A} \in \Alg{\omega}$ with fine spectrum
       $\fS{{\bf A}}$ that has the property
       \[
        \exists n \in \N, n \ge 2 : \ \aS{{\bf A}}{n} = 1 \ \textrm{ and } \ \exists m \in \N, m > n : \ \aS{{\bf A}}{m} > 1.
       \]
       Then define the sequence $\left( \Sigma_i \right)_{i \in \N}$ in $\L{p}$ as follows:
       \[
        \Sigma_i := \begin{cases}
                     \B{i}{p} \times \B{i}{p}, & \textrm{if } i < n; \\
                     \fSn{{\bf A}}{i}, & \textrm{otherwise}.
                    \end{cases}
       \]
       It is easy to see that the corresponding $\Sigma \subseteq \L{p}$ is an equivalence relation
       that is closed under the implication operator. Therefore, there exists a groupoid ${\bf B} \in \Alg{\omega}$ with
       $\fS{{\bf B}} = \Sigma$. This is a contradiction to (\ref{GenAssLawB}) because $\aS{{\bf B}}{2} = 1$ but $\aS{{\bf B}}{m} > 1$.
\end{enumerate}
\eproof

\section{Equivalent representation of bracketings}

In this section we are going to define an equivalent representation of bracketings.
First let us look at the number of bracketings with a given occurence number.
\bd\label{DefCatalan}
The {\it generalized {\sc Catalan} numbers} $\C{n}{p}$ ($n \ge 0, p \ge 1$) are defined by the following recursion:
\begin{itemize}
 \item $\C{0}{p} := 1$;
 \item $\displaystyle \C{n}{p} := \sum_{i_1, \ldots, i_p \in \N, \sum_{k=1}^{p} i_k = n - 1}
        \left( \prod_{l=1}^{p} \C{i_l}{p} \right)$.
\end{itemize}

With respect to the definition of the bracketings as terms of $\Tp$ and the properties of the occurence number stated in Section
\ref{sectionOne} it follows:
\[
 \forall n \in \N : \ \left| \B{n}{p} \right| = \C{n}{p}.
\]
We know from \cite{Klarner} that
\[
 \C{n}{p} = \frac{1}{(p-1) \cdot n + 1} \cdot \binom{p \cdot n}{n}
\]
holds for all $n \in \N$. We will prove this in a more general form in Theorem~\ref{ThmMNKP}.
\ed
With the same idea as above and the fact that $\fS{{\bf A}}$ is a congruence relation in
$\Tp$ we see that for any groupoid ${\bf A}$:
\[
 \forall n \in \Nplus : \ \aS{{\bf A}}{n} \le \sum_{i_1, \ldots, i_p \in \N, \sum_{k=1}^{p} i_k = n - 1}
        \left( \prod_{l=1}^{p} \aS{{\bf A}}{i_l} \right).
\]
Now we are going to define an equivalent representation for bracketings which will be useful to prove the explicit formula of the
generalized {\sc Catalan} numbers.
\bd\label{DefInsertionTuple}
The {\it insertion tuple} of a bracketing $t$ is the tuple $\ST \left( t \right) \in \N^{\occN{t}}$ whose
$i$-th entry is one plus the number of $x$'s preceding the $i$-th symbol $\omega$ in the (prefix notation of) $t$.
It can be also defined recursively as follows.
\begin{itemize}
 \item $\ST \left( x \right) := \emptyset \in \N^{0}$;
 \item $\ST \left( \omega t_1 \ldots t_p \right) := \left( 1, {\bf v}^1, \ldots, {\bf v}^p \right)$ 
       is the consecutive sequence of $1$ and the ${\bf v}^i$ where ${\bf v}^i$ is the insertion tuple
       of $t_i$ with an additional shift that is added to each entry of the tuple:
       \[ {\bf v}^i := \ST \left( t_i \right) + \sum_{k=1}^{i-1} \length{t_k} \in \N^{\occN{t_i}}.\]
       The shift is exactly the sum of the lengths of the previous bracketings just as in Lemma~\ref{LemmaBeta}.
\end{itemize}
\ed
The insertion tuples of the first binary bracketings have been presented in Section~\ref{sectionOne}.
We introduce the following notation, which we will need later.
For any $k \in \Nplus$ and $n \in \N$ let
\[ \M{n}{k}{p} :=
        \rset{{\bf u} \in \N^n}{
         1 \le u_i \le (p-1) \cdot (i-1) + k
         \textrm{ and } u_i \le u_j \ \ (1 \le i \le j \le n)}.\]
The insertion tuples can be characterized as follows.
\bpro\label{PropST}
For $t \in \B{n}{p}$ the following hold:
\begin{enumerate}[(a)]
 \item $\ST \left( t \right) \in \N^{n}$.
 \item\label{IT_Part_b} If we know the insertion tuple ${\bf u} := \ST \left( t \right)$
       of $t$, then we also know the insertion tuple ${\bf v} := \ST \left( \beta_i \left( t \right) \right)$
       of $\beta_i \left( t \right)$ (for $i = 1, \ldots, \l{p} (n)$):
       \[
        v_q =
        \begin{cases}
         u_q, & \textrm{if } q \in \set{1, \ldots, l}; \\
         i,   & \textrm{if } q = l + 1; \\
         u_{q-1} + p-1, & \textrm{if } q \in \set{l+2, \ldots, n+1},
        \end{cases}
       \]
       where $l := \max \set{0} \cup \lset{q \in \set{1, \ldots, n}}{u_q \le i}$.
 \item $\ST \left[ \B{n}{p} \right] = \M{n}{1}{p}$
 \item $\ST$ is an injective map.
 \item\label{IT_Part_e} The name insertion tuple is justified, because with ${\bf u} := \ST \left( t \right)$
       and the two definitions $t_0 := x$, $t_i := \beta_{u_i} \left( t_{i-1} \right)$ (for $i \in \set{1, \ldots, n}$)
       we have $t = t_n$.
\end{enumerate}
\epro
\bproof
\begin{enumerate}[(a)]
 \item Clear from the definition.
 \item Remember that $\beta_i (t)$ is the insertion of $\omega x \ldots x$ at the $i$-th symbol $x$  in $t$ and
       that $u_q - 1$ equals the number of $x$'s preceding the $q$-th symbol $\omega$.
       The statement becomes clear if we observe that $l$ is the position of the last symbol $\omega$ having less than $i$ many
       $x$'s before it.
 \item \begin{itemize}
        \item[``$\subseteq$'':] Let $t \in \B{n}{p}$ be an arbitrary bracketing and let ${\bf u} := \ST (t)$.
                                It is clear from the definition that $1 \le u_i$ and ${\bf u}$ is monotone. For the upper bound
                                let us consider the $i$-th occurence of the symbol $\omega$ in $t$. If we delete all symbols from this
                                $\omega$ to the end of the bracketing then we obtain the prefix of a bracketing with occurence number
                                $i-1$ with at least one $x$ missing. Therefore, the number of $x$'s preceding this $\omega$ is at most
                                $\l{p} (i-1) - 1$. Hence, $u_i \le \l{p} (i-1)$.
        \item[``$\supseteq$'':] We use induction on $n$ to show this. The case $n=0$ is clear, because both sides are $\set{\emptyset}$.
                                For the induction step from $n$ to $n+1$ let ${\bf u} \in \M{n+1}{1}{p}$. Then
                                $(u_q)_{q=1,\ldots,n} \in \M{n}{1}{p}$. By induction hypothesis we have a bracketing
                                $t \in \B{n}{p}$ such that $\ST \left( t \right) = (u_q)_{q=1,\ldots,n}$. With the help of
                                (\ref{IT_Part_b}) and the monotonicity of $\ST \left( t \right)$ we see that
                                \[
                                 \ST \left( \beta_{u_{n+1}} \left( t \right) \right) = {\bf u}.
                                \]
          \end{itemize}
 \item The entries of the insertion tuples tell us the positions of the symbols $\omega$.
 \item It can be shown by induction on $i$ with the help of (\ref{IT_Part_b}) that
       \[
        \ST ( t_i ) = (u_q)_{q=1,\ldots,i} \quad (i = 1, \ldots, n)
       \]
       holds. So we have $\ST (t_n) = {\bf u} = \ST (t)$. Then by injectivity $t = t_n$ follows.
\end{enumerate}
\eproof
In \cite{HiltonPederson} a bijection is given between bracketings and certain lattice paths called $p$-good paths.
We invite the reader to find a bijection between $p$-good paths and insertion tuples.
The sets $\M{n}{k}{p}$ generalize insertion tuples, and the corresponding paths generalize the
$p$-good paths by shifting the bounding line $k-1$ steps up. Therefore, Theorem~\ref{ThmMNKP} can be considered as a generalization of
Theorem 0.4 from \cite{HiltonPederson}.

\bl\label{LemmaMNKP}
For any $k \in \Nplus$ and $n \in \N$ we have:
\begin{enumerate}[(a)]
 \item $\left| \M{0}{k}{p} \right| = 1$;
 \item $\displaystyle \left| \M{n+1}{k}{p} \right| = \sum_{l=0}^{k-1} \left| \M{n}{p+l}{p} \right|$.
\end{enumerate}
\el
\bproof
(a) is trivial. For (b) we partition the set $\M{n+1}{k}{p}$ into the disjoint subsets
\[
 S_l := \lset{{\bf u} \in \M{n+1}{k}{p}}{u_1 = l} \quad (l = 1, \ldots, k).
\]
It is easy to verify that the map
\[
 \function{\varphi_l}{S_l}{\M{n}{k+p-l}{p}}{{\bf u}}{(u_q - (l - 1))_{q=2,\ldots,n+1}}
\]
is a bijection for $l = 1, \ldots, k$. Therefore

\noindent
\parbox{\textwidth}{\[
  \left| \M{n+1}{k}{p} \right| = \sum_{l=1}^{k} \left| \M{n}{k+p-l}{p} \right| = \sum_{l=0}^{k-1} \left| \M{n}{p+l}{p} \right| .
\]
\eproof}
The following lemma can be shown by a straightforward induction on $k$.
\bl\label{LemmaProduct}
For any $k \in \Nplus$ and $n \in \N$ we have:
\[
 \begin{array}{rl}
  & \displaystyle \sum_{m=0}^{k-1} \frac{p + m}{(p - 1) \cdot (n + 1) + m + 1} \cdot \prod_{l=m+1}^{n+m} \left( (p - 1) \cdot (n+1) + l \right) \\
  = & \displaystyle \frac{1}{n+1} \cdot \frac{k}{(p - 1) \cdot (n+1) + k} \cdot \prod_{l=k}^{n+k} \left( (p - 1) \cdot (n+1) + l \right).
 \end{array}
\]
\el

\bthm\label{ThmMNKP}
For any $k \in \Nplus$ and $n \in \N$ we have:
\[
 \begin{array}{rcl}
  \left| \M{n}{k}{p} \right| & = & \displaystyle
    \frac{1}{n!} \cdot \frac{k}{(p - 1) \cdot n + k} \cdot \prod_{l=k}^{n+k-1} \left( (p - 1) \cdot n + l \right) \\
  & = & \displaystyle \frac{k}{(p-1) \cdot n + k} \cdot \binom{p \cdot n + k - 1}{n}.
 \end{array}
\]
\ethm
\bproof
It is a routine induction proof using the recursion formula from
Lemma~\ref{LemmaMNKP} and the previous Lemma~\ref{LemmaProduct}.
\eproof

\section{Polynomial spectra}\label{sectionFour}

In this section we give three different examples of polynomial spectra which solve problem 3 in \cite{CsakanyWaldhauser}.

\bex\label{ExPolyTamas}
Let $k$ be a fixed natural number. We define an equivalence relation
$\Sigma \subseteq \L{2}$ as follows. For a bracketing
$s=\omega t_{1}t_{2}\in \B{n}{2}$ we call $t_{1}$ the
\emph{left factor of }$s$ and denote it by $\operatorname{left}\left(
s\right)  $, and we put $\operatorname{left}\left(  x\right)  =x$ (cf. \cite{CsakanyWaldhauser}).
For $\eq{s}{t} \in \B{n}{2} \times \B{n}{2}$ let
\[
\eq{s}{t} \in \Sigma :\Longleftrightarrow \left\vert \operatorname{left}^{i}\left(
s\right)  \right\vert =\left\vert \operatorname{left}^{i}\left(  t\right)
\right\vert \text{ for }i=1,\ldots,k.
\]
The set $\Sigma$ is closed under the implication operator, thus it appears as
the fine spectrum of some groupoid ${\bf A}$. The corresponding associative spectrum is
a polynomial of degree $k$:
\[
\aS{{\bf A}}{n}=\binom{n-1}{k}+\binom{n-1}{k-1}+\cdots+\binom{n-1}{1}+\binom{n-1}{0}.
\]
\eex
\bproof
It is straightforward to check that $\delta_n \left(  \Sigma_{n}\right)
\subseteq\Sigma_{n+1}$ holds for all $n \in \N$. Let $s\in
\B{n}{2}$ be an arbitrary bracketing, and let us abbreviate
$\left\vert \operatorname{left}^{i}\left(  s\right)  \right\vert $ by $l_{i}$.
Clearly we have%
\[
1\leq l_{k}\leq l_{k-1}\leq\ldots\leq l_{2}\leq l_{1}\leq n,
\]
where the inequalities are strict, except maybe for a couple of repeated $1$s
at the beginning. We have to count how many such $k$-tuples exist. If the
number of $1$s at the beginning is $i$, then we have to choose $k-i$ different
numbers from the set $\left\{  2,\ldots,n\right\}  $, hence the number of
possibilities is $\binom{n-1}{k-i}$. Thus we have%
\[
\left\vert \B{n}{2} / \Sigma_{n} \right\vert =\binom{n-1}
{k}+\binom{n-1}{k-1}+\cdots+\binom{n-1}{1}+\binom{n-1}{0}\text{,}
\]
which is indeed a polynomial of degree $k$.
\eproof

\bex
Let $k \in \Nplus$ be an integer and define the relation $\Sigma \subseteq \L{p}$ as follows: For an identity
$\eq{s}{t} \in \B{n}{p} \times \B{n}{p}$ denote by ${\bf u} := \ST (s)$ and ${\bf v} := \ST (t)$ the
insertion tuples of the bracketings. We define $\Sigma$ by
\[
 \eq{s}{t} \in \Sigma : \Longleftrightarrow \begin{cases}
                                                   s = t, & \textrm{if } n < k; \\
                                                   \forall i \in \set{n - k + 1, \ldots, n} : \ u_i = v_i, & \textrm{if } n \ge k.
                                                  \end{cases}
\]
It holds that $\Sigma$ is an equivalence relation that is closed under the implication operator. Therefore, there exists
a groupoid ${\bf A}$ such that $\fS{{\bf A}} = \Sigma$. The associative spectrum of ${\bf A}$ is
\[
 \aS{{\bf A}}{n} = \begin{cases}
                    \C{n}{p}, & \textrm{if } n < k; \\
                    \displaystyle \frac{(p-1) \cdot (n-k) + 1}{(p-1) \cdot n + 1} \cdot \binom{(p-1) \cdot n + k}{k}, & \textrm{if } n \ge k.
                   \end{cases}
\]
\eex
\bproof
Remember that  $\eq{s}{t} \in \Sigma$ means that the last $k$ entries of the insertion tuples of $s$ and $t$ are equal, or equivalently
that the last $k$ symbols $\omega$ are in the same place in $s$ and $t$. Therefore,
it is easy to verify that $\delta_n \left( \Sigma_{n} \right) \subseteq \Sigma_{n+1}$ holds for all $n \in \N$.

%
%
%

To know the associative spectrum we have to count the insertion tuples with different last $k$ entries. From Proposition
\ref{PropST} we know that $\ST \left[ \B{n}{p} \right] = \M{n}{1}{p}$. If we look at the last $k \le n$ entries we see that
they form exactly the set $\M{k}{(p-1) \cdot (n-k) + 1}{p}$. So the formula can be obtained from Theorem~\ref{ThmMNKP}.
\eproof

\bex
The binary groupoid ${\bf G}:= \pair{\Z_6[Y]}{\oplus}$ (where $\Z_6[Y]$ is the polynomial ring over $\Z_6$ in the variable $Y$)
with the binary operation
\[
 \function{\oplus}{\Z_6[Y]^2}{\Z_6[Y]}{(X_1, X_2)}{3Y \cdot X_1 + 2Y \cdot X_2}
\]
has the associative spectrum
\[
 \forall n \in \N, n \ge 2 : \ \aS{{\bf G}}{n} = \frac{n^2 + n - 2}{2} .
\]
Instead of $\Z_6[Y]$ another ring can be chosen which has zero divisors (in this case $3Y$ and $2Y$) whose
powers are all different.
\eex
\bproof
In \cite{CsakanyWaldhauser} the notion of left and right depth sequence were defined. Here, we only need two special cases: for
$s \in \B{n}{2}$ let $d_l (s)$ denote the left depth of the first variable of $s$, and
let $d_r (s)$ denote the right depth of the last variable of $s$.
On the binary tree corresponding to $s$ one can see $d_l (s)$ as the length of the path connecting the root and the leftmost
leaf and $d_r (s)$ as the length of the path connecting the root and the rightmost leaf.
From this interpretation it is clear that for all $t_1, t_2 \in \T$:
\begin{equation}\label{SomeEquation}
 d_l (\omega t_1 t_2) = d_l (t_1) + 1 \ \textrm{ and } \ d_r (\omega t_1 t_2) = d_r (t_2) + 1 .
\end{equation}
Later it will be useful to compute these numbers from the insertion tuple ${\bf u} := \ST (s)$:
\begin{align*}
 d_l (s) & = \left| \Big\{ q \in \set{1, \ldots, n} \ \Big| \ u_q = 1 \Big\} \right| , \\
 d_r (s) & = \left| \rset{q \in \set{1, \ldots, n}}{u_q = \l{2} (q-1)} \right| .
\end{align*}
It is a routine induction using (\ref{SomeEquation}) to check that for $n > 0$
\[
 \left( \en{1} (s) \right)^{{\bf G}} (X_1, \ldots, X_{n+1}) = (3Y)^{d_l (s)} \cdot X_1 + (2Y)^{d_r (s)} \cdot X_{n+1}
\]
holds. Remember that the length function is $\l{2} (n) = n+1$ and that $\en{1} (s)$ contains exactly the variable symbols
$x_1, \ldots, x_{\l{2} (n)}$ such that $\en{1} (s)$ can be interpreted as a ($n+1)$-ary term.

From this it follows that the fine spectrum of ${\bf G}$ can be characterized
as:
\[
 \eq{s}{t} \in \fS{{\bf G}} \Longleftrightarrow  d_l (s) =  d_l (t) \textrm{ and } d_r (s) =  d_r (t) .
\]
To gain the formula for the associative spectrum we have to count all possibilities for $d_l (s)$ and $d_r (s)$.
For $n \ge 2$ we have the following restrictions:
\[
 d_l (s) \ge 1 \textrm{, } d_r (s) \ge 1 \textrm{ and } 3 \le d_l (s) + d_l (r) \le n + 1.
\]
This is pretty clear using the insertion tuple because the first entry $u_1 = 1$ counts for both depths, the second entry $u_2$ can either
be 1 or $2 = \l{2} (1)$ and the other entries $u_k$ ($k = 3, \ldots, n$) can be 1, $\l{2} (k - 1)$ or something in between.
It is also clear that all such possibilities can occur. So we have
\[
 \left( \sum_{l=1}^{n} n + 1 - l \right) - 1 = \left( \sum_{k=1}^{n} k \right) - 1 = \frac{n^2 + n - 2}{2}
\]
possibilities.
\eproof

\section{The lattice of fine spectra}

The {\sc Galois}-closed sets $\IdBrack \K$ ($\K \subseteq \Alg{\omega}$) form a complete lattice as a closure system in $\p{\L{p}}$.
We will denote this complete lattice by $\FSB = \left(\FS, \bigwedge, \bigvee \right)$ (FS stands for fine spectra
because we know from Corollary~\ref{TheoremCharFS} that the elements of this lattice are exactly all fine spectra).
To unterstand associative spectra it would be very useful to unterstand this lattice. As a beginning we will
look at the covering relation $\prec$.

\bpro
For any $\fS{{\bf A}}, \fS{{\bf B}} \in \FS$ the following holds:
\[
 \begin{array}{rcl}
  \fS{{\bf A}} \prec \fS{{\bf B}} & \Longleftrightarrow & \exists ! n \in \N : \ \fSn{{\bf A}}{n} \ne \fSn{{\bf B}}{n}
  \textrm{, and for this $n$ we have } \\
  & & \fSn{{\bf A}}{n} \prec \fSn{{\bf B}}{n} \textrm{ in the lattice } \Eq \left( \B{n}{p} \right).
 \end{array}
\]
\epro
\bproof
It is clear that the condition on the right is sufficient. For the necessity let us assume that $\fS{{\bf A}} \prec \fS{{\bf B}}$.
If $\fS{{\bf A}}$ and $\fS{{\bf B}}$ differ at least at two positions, say $\fSn{{\bf A}}{n} \ne \fSn{{\bf B}}{n}$ and
$\fSn{{\bf A}}{m} \ne \fSn{{\bf B}}{m}$ for some $n < m \in \N$, then
\[
 \Sigma_k := \begin{cases}
              \fSn{{\bf A}}{k}, & \textrm{if } k \le n; \\
              \fSn{{\bf B}}{k}, & \textrm{if } k > n
             \end{cases}
\]
defines a fine spectrum that is strictly between $\fS{{\bf A}}$ and $\fS{{\bf B}}$ contradicting that $\fS{{\bf A}} \prec \fS{{\bf B}}$.
If $\fS{{\bf A}}$ and $\fS{{\bf B}}$ differ only at one position, say $\fSn{{\bf A}}{n} \ne \fSn{{\bf B}}{n}$ and
$\fSn{{\bf A}}{n} \nprec \fSn{{\bf B}}{n} \textrm{ in the lattice } \Eq \left( \B{n}{p} \right)$, then
\[
 \Sigma_k := \begin{cases}
              \fSn{{\bf A}}{k}, & \textrm{if } k \ne n; \\
              \pi, & \textrm{if } k = n
             \end{cases}
\]
defines a fine spectrum that is strictly between $\fS{{\bf A}}$ and $\fS{{\bf B}}$ if $\pi \in \Eq \left( \B{n}{p} \right)$ is an
equivalence relation such that $\fSn{{\bf A}}{n} < \pi < \fSn{{\bf B}}{n}$.
\eproof
As a consequence we the obtain the following characterization of the atoms and coatoms of $\FSB$:
\bcor
There are no atoms in $\FSB$. For any $\fS{{\bf A}} \in \FS$ we have:
\[
  \fS{{\bf A}} \textrm{ is a coatom in } \FSB \Longleftrightarrow \forall n \in \N \setminus \set{2} : \ \aS{{\bf A}}{n} = 1 \ \textrm{ and } \
  \aS{{\bf A}}{2} = 2.
\]
\ecor
\bproof
%
The description for the coatoms follows from the above proposition.
For the atoms we assume that $\fS{{\bf A}} \in \FS$ is an atom in $\FSB$. From the previous proposition we see that
$\fSn{{\bf A}}{n}$ is the equality relation on $\B{n}{p}$ for all but one $n \in \N$. It is clear that such a $\fS{{\bf A}}$ cannot
be closed under the implication operator.
\eproof
From the previous corollary we know that the number of coatoms is exactly the number of possibilities to group
$p$ elements into two classes (because $\left| \B{2}{p} \right| = p$). So we get:
\bcor
There are exactly $2^{p-1} - 1$ coatoms in $\FSB$.
\ecor
We prove that the cardinality of the set of sequences of natural numbers that
arise as associative spectra is continuum. Clearly, it cannot be more, so it
suffices to construct continuously many different spectra, and it suffices to
do it in the binary case. First we need a definition: if $\omega xx=\left(
xx\right)  $ is a subbracketing of $s\in \B{n}{2}$, then we
say that $\left(  xx\right)  $ is a \emph{pair of eggs} in $s$. (Actually the
two $x$'s are the eggs, see \cite{CsakanyWaldhauser}.)

\bl
Let $\tau_{n}$ be the equivalence relation on $\B{n}{2}$,
where the bracketings with at least $3$ pairs of eggs form one class, and all
the other bracketings are singletons. Then $\tau_{n} \supsetneq \delta_{n-1}\left(  \tau
_{n-1}\right)  $ for all $n\geq5$.
\el
\bproof
The operators $\gamma_{i},\beta_{i}$ do not decrease the number of eggs, hence
$\delta_{n-1}\left(  \tau_{n-1}\right)  \subseteq\tau_{n}$ for all $n\geq1$. For
every $n\geq5$ one can find a bracketing with occurence number $n$, which
cannot be obtained by these operators from any bracketing (with occurence
number $n-1$) with at least three pairs of eggs. For example,
\[
t=\left(  \cdots\left(  \left(  \left(  xx\right)  \left(  xx\right)  \right)
x\right)  x\cdots x\right)  \left(  xx\right)
\]
\parbox{\textwidth}{
is such a bracketing. Thus $t$ is a singleton in $\delta_{n-1}\left(  \tau
_{n-1}\right)  $, but not a singleton in $\tau_{n}$. This shows that
$\delta_{n-1}\left(  \tau_{n-1}\right)  \neq\tau_{n}$ if $n\geq5$.
\eproof}

\bthm\label{ThmUncount}
There exist $2^{\aleph_{0}}$ different associative spectra.
\ethm
\bproof
Let $S$ be the set of 0--1 sequences whose first five entries are 0. For every
$\mathbf{a}=\left\{  a_{n}\right\}  _{n=0}^{\infty}\in S$ we construct a
sequence of equivalence relations $\sigma_{n}^{\mathbf{a}}\subseteq
\B{n}{2}\times \B{n}{2}$ recursively:
\[
\sigma_{n}^{\mathbf{a}}=\left\{
\begin{array}
[c]{cc}%
\delta_{n-1}\left(  \sigma_{n-1}^{\mathbf{a}}\right),  & \text{if }a_{n}=0;\\
\tau_{n}, & \text{if }a_{n}=1.
\end{array}
\right.
\]
Note that we do not have to define the ``initial
value'' $\sigma_{0}^{\mathbf{a}}$ since $\B{0}{2}$ is a one-element set. Observe also that $\sigma_{n}^{\mathbf{a}}$
is the equality relation on $\B{n}{2}$ for $n\leq4$ for
every $\mathbf{a}\in S$.

First we claim that $\sigma_{n}^{\mathbf{a}}\subseteq\tau_{n}$ for every
$\mathbf{a}\in S$. This is clear for $n=0$ (and also for $n=1,2,3,4$), and
then we can proceed by induction. Suppose that $\sigma_{n-1}^{\mathbf{a}%
}\subseteq\tau_{n-1}$. If $a_{n}=1$, then $\sigma_{n}^{\mathbf{a}}=\tau_{n}$;
if $a_{n}=0$, then $\sigma_{n}^{\mathbf{a}}=\delta_{n-1}\left(  \sigma
_{n-1}^{\mathbf{a}}\right)  \subseteq \delta_{n-1}\left(  \tau_{n-1}\right)
\subsetneq \tau_{n}$ by the previous lemma, and by the monotonicity of $\delta_{n-1}$.

Now we can verify that $\sigma^{\mathbf{a}}$ is a fine spectrum: if $a_{n}=0$,
then $\sigma_{n}^{\mathbf{a}}=\delta_{n-1}\left(  \sigma_{n-1}^{\mathbf{a}}\right)
$; if $a_{n}=1$, then $\sigma_{n}^{\mathbf{a}}=\tau_{n}\supsetneq\delta_{n-1}\left(
\tau_{n-1}\right)  \supseteq\delta_{n-1}\left(  \sigma_{n-1}^{\mathbf{a}}\right)  $,
hence Corollary~\ref{TheoremCharFS} applies.

We need to check yet that different elements of $S$ give different associative
spectra. Let $\mathbf{a}\neq\mathbf{b}\in S$, and suppose that $a_{i}\neq
b_{i}$, say $a_{i}=0$ and $b_{i}=1$. Then we have $\sigma_{i}^{\mathbf{a}%
}=\delta_{i-1}\left(  \sigma_{i-1}^{\mathbf{a}}\right)  \subseteq\delta_{i-1}\left(
\tau_{i-1}\right)  \subsetneq\tau_{i}=\sigma_{i}^{\mathbf{b}}$. We have proved
that $\sigma_{i}^{\mathbf{a}}\subsetneq\sigma_{i}^{\mathbf{b}}$, and this means
that not only the two fine spectra, but the corresponding spectra are also
different: $\left\vert \B{i}{2}/\sigma_{i}^{\mathbf{a}%
}\right\vert >\left\vert \B{i}{2}/\sigma_{i}^{\mathbf{b}%
}\right\vert $.
\eproof

\brem
From the previous proof we see that
\[
 \forall {\bf a}, {\bf b} \in S : \ \sigma^{\mathbf{a}} \subseteq \sigma^{\mathbf{b}} \Longleftrightarrow {\bf a} \le {\bf b},
\]
which means that $S$ embeds into $\FSB$ as a poset (not as a lattice!). Since $S$ is isomorphic to $\p{\N}$, we have
$\p{\N}$ as a subposet in $\FSB$.
On the other hand clearly $\FSB$ embeds into $\p{\L{p}} \cong \p{\N}$. Therefore, $\FSB$ and $\p{\N}$ are {\it equimorphic}.
This shows for example, that there is a chain and an antichain of continuum cardinality in $\FSB$.
\erem

\section{Spectra of finite groupoids}

There are only countably many finite groupoids, hence Theorem \ref{ThmUncount} shows
that there are spectra which can be realized only on infinite groupoids. It
would be interesting to see, under what conditions a (fine) spectrum
is realizable on a finite groupoid. One obvious necessary condition: the spectrum
has to be recursive (computable by a Turing machine). If there exists
$N\in\N$ such that $\sigma_{n}=\delta_{n-1}\left(  \sigma
_{n-1}\right)  $ holds for all $n>N$, then the sequence $\sigma_{n}$ is
recursive. We conjecture that this condition is sufficient in order to realize
a fine spectrum on a finite groupoid (cf. Proposition \ref{PropFinallyAssoc}).
The condition is not necessary, as the following example shows.

\bex
We construct a finite groupoid with the \textquotedblleft three-egg
spectrum\textquotedblright\ $\tau_{n}$. First let us consider the groupoid
$\mathbf{A}$ given by the following multiplication table.%
\[%
\begin{tabular}
[c]{c|cccc}
& $0$ & $1$ & $2$ & $3$\\\hline
$0$ & $0$ & $0$ & $0$ & $0$\\
$1$ & $0$ & $0$ & $0$ & $1$\\
$2$ & $0$ & $0$ & $1$ & $2$\\
$3$ & $0$ & $1$ & $2$ & $2$%
\end{tabular}
\ \ \ \
\]
One can prove by induction that for any bracketing $s$ the maximal value of
the corresponding term function $s^{\mathbf{A}}$ is $\max\left(  3-e,0\right)
$ where $e$ is the number of pairs of eggs in $s$. This maximal value is
attained for example at $s^{\mathbf{A}}\left(  3,\ldots,3\right)  $. This
shows that $\sigma_{n}\left(  \mathbf{A}\right)  \supseteq\tau_{n}$, since
bracketings with at least three pairs of eggs induce constant $0$ term
functions. (Actually one can verify that $\mathbf{A}\models s_{1}\approx
s_{2}$ iff either both $s_{1}$ and $s_{2}$ contain at least three pairs of
eggs, or both contain at most two pairs of eggs and these are at the same
positions in $s_{1}$ and $s_{2}$.)

In order to isolate bracketings with at most two eggs, we blow up the nonzero
elements of $\mathbf{A}$ using the Sheffer operation on the two-element set.
We present this operation with somewhat unusual notation:%
\[
\renewcommand{\arraystretch}{1.5}%
\begin{tabular}
[c]{r|rr}
& $\widehat{\square}$ & $\widetilde{\square}$\\\hline
$\widehat{\square}$ & $\widetilde{\square}$ & $\widehat{\square}$\\
$\widetilde{\square}$ & $\widehat{\square}$ & $\widehat{\square}$%
\end{tabular}
\ \
\]
We replace each nonzero element of $\mathbf{A}$ with two elements: one wearing
a hat, the other one wearing a tilde, and we define the multiplication such
that the numbers get multiplied as in $\mathbf{A}$, and headgears get
multiplied according to the Sheffer operation. We obtain the following
seven-element groupoid $\widehat{\mathbf{A}}$:%
\[
\renewcommand{\arraystretch}{1.6}\begin{tabular}
[c]{p{0.3cm}|p{0.3cm}:p{0.3cm}p{0.3cm}:p{0.3cm}p{0.3cm}:p{0.3cm}p{0.3cm}}
& $0$ & $\widehat{1}$ & $\widetilde{1}$ & $\widehat{2}$ & $\widetilde{2}$ &
$\widehat{3}$ & $\widetilde{3}$\\\hline
$0$ & $0$ & $0$ & $0$ & $0$ & $0$ & $0$ & $0$\\\hdashline
$\widehat{1}$ & $0$ & $0$ & $0$ & $0$ & $0$ & $\widetilde{1}$ & $\widehat
{1}$\\
$\widetilde{1}$ & $0$ & $0$ & $0$ & $0$ & $0$ & $\widehat{1}$ & $\widehat
{1}$\\\hdashline
$\widehat{2}$ & $0$ & $0$ & $0$ & $\widetilde{1}$ & $\widehat{1}$ &
$\widetilde{2}$ & $\widehat{2}$\\
$\widetilde{2}$ & $0$ & $0$ & $0$ & $\widehat{1}$ & $\widehat{1}$ &
$\widehat{2}$ & $\widehat{2}$\\\hdashline
$\widehat{3}$ & $0$ & $\widetilde{1}$ & $\widehat{1}$ & $\widetilde{2}$ &
$\widehat{2}$ & $\widetilde{2}$ & $\widehat{2}$\\
$\widetilde{3}$ & $0$ & $\widehat{1}$ & $\widehat{1}$ & $\widehat{2}$ &
$\widehat{2}$ & $\widehat{2}$ & $\widehat{2}$%
\end{tabular}
\]
We did not blow up $0$, hence bracketings with at least three pairs of eggs
still induce constant term functions, and thus we have $\sigma_{n}%
\bigl(%
\widehat{\mathbf{A}}%
\bigr)%
\supseteq\tau_{n}$. On the other hand, if $s$ contains at most two pairs of
eggs, then $s^{\mathbf{A}}\left(  x_{1},\ldots,x_{n}\right)  \neq0$ if
$x_{1},\ldots,x_{n}\in\big\{\widehat{3},\widetilde{3}\big\}$. This means that
substituting $\widehat{3}$s and $\widetilde{3}$s into $s^{\mathbf{A}}$ we can
recover all information about hats and tildes, that is we can determine the
term function corresponding to $s$ over the Sheffer operation. This operation
is Catalan, hence from the term function we can recover the bracketing.
Consequently $s$ is a singleton in $\sigma_{n}%
\bigl(%
\widehat{\mathbf{A}}%
\bigr)%
$, hence $\sigma_{n}%
\bigl(%
\widehat{\mathbf{A}}%
\bigr)%
=\tau_{n}$.
\eex
In Section \ref{sectionFour} we gave examples for polynomial spectra using Corollary \ref{TheoremCharFS}.
The groupoids that we obtained this way were infinite
groupoids of the form $\Tp / \Sigma$, but below we will construct
a finite groupoid with a polynomial spectrum.
\bex
Let us define a binary operation on the set $A=\left\{  0,1,\ldots
,k+1\right\}  $ by
\[
x\cdot y=\left\{  \!\!\!%
\begin{tabular}
[c]{ll}%
$0,$ & if $x=0;$\\
$1,$ & if $x\neq0=y;$\\
$\min\left(  x+1,k+1\right)  ,$ & if $x\neq0\neq y.$%
\end{tabular}
\ \right.
\]
The associative spectrum of the groupoid $\mathbf{A}\mathbb{=}\left(
A;\cdot\right)  $ is a polynomial of degree $k$.
\eex
\bproof
Let $\Sigma$ be the equivalence relation defined in Example \ref{ExPolyTamas}. We prove that
the fine spectrum of $\mathbf{A}$ is $\sigma_{n}\left(  \mathbf{A}\right)
=\Sigma_{n}$. To avoid notational difficulties we prove it only for $k=3$; it
will be clear from the proof how the construction works for arbitrary $k$. To
have a better view of the operation, let us write out the multiplication
table.%
\[%
\begin{tabular}
[c]{r|rrrrr}
& $0$ & $1$ & $2$ & $3$ & $4$\\\hline
$0$ & $0$ & $0$ & $0$ & $0$ & $0$\\
$1$ & $1$ & $2$ & $2$ & $2$ & $2$\\
$2$ & $1$ & $3$ & $3$ & $3$ & $3$\\
$3$ & $1$ & $4$ & $4$ & $4$ & $4$\\
$4$ & $1$ & $4$ & $4$ & $4$ & $4$%
\end{tabular}
\
\]
First we prove that $\sigma_{n}\left(  \mathbf{A}\right)  \subseteq\Sigma_{n}$
for all $n\in\mathbb{N}\hspace{0cm}$. It suffices to show that for any $t\in
B_{n}^{\left(  2\right)  }$, the values of $l_{1},l_{2},l_{3}$ can be read off
from the term function $t^{\mathbf{A}}$ corresponding to $t$. Taking a look at
the multiplication table, we see immediately that $x\cdot y=0$ iff $x=0$.
Therefore a product of arbitrarily many elements (with arbitrarily inserted
parentheses) equals $0$ iff the first (i.e. leftmost) element is $0$. We
record this fact with the following (hopefully intuitive) notation, where
$\ast$ symbolises an arbitrary nonzero element:%
\begin{equation}\label{nulla}
\left(  0\cdots\right)  =0,\quad\left(  \ast\cdots\right)  \neq0.
\end{equation}
From this observation and from the idempotence of the element $4$ we infer%
\begin{align*}
\left(  4\cdots4\right)  \cdot\left(  4\cdots0\right)   &  =4\cdot\ast=4,\\
\left(  4\cdots4\right)  \cdot\left(  0\cdots0\right)   &  =4\cdot0=1.
\end{align*}
This means that $l_{1}$ can be computed from the values of $t^{\mathbf{A}}$:%
\[
l_{1}=\max\Big\{i~\Big|~t^{\mathbf{A}}(\underset{i}{\underbrace{4,\ldots,4}%
},0,\ldots0)=1\Big\}.
\]
Knowing the value of $l_{1}$, we can find $l_{2}$ using the following
observations:%
\begin{align*}
\left(  \left(  4\cdots4\right)  \cdot\left(  4\cdots0\right)  \right)
\cdot\left(  4\cdots4\right)   &  =\left(  4\cdot\ast\right)  \cdot
4=4\cdot4=4,\\
\left(  \left(  4\cdots4\right)  \cdot\left(  0\cdots0\right)  \right)
\cdot\left(  4\cdots4\right)   &  =\left(  4\cdot0\right)  \cdot4=1\cdot4=2.
\end{align*}
Hence $l_{2}$ can be recovered from $t^{\mathbf{A}}$ as%
\[
l_{2}=\max\Big\{i~\Big|~t^{\mathbf{A}}%
\bigl(%
\overset{l_{1}}{\overbrace{\underset{i}{\underbrace{4,\ldots,4}},0,\ldots0}%
},4,\ldots,4%
\bigr)%
=2\Big\}.
\]
Note that if $l_{1}=1$, then we cannot make such substitutions, but in this
case clearly $l_{2}=l_{3}=1.$

Similarly, $l_{3}$ can be obtained, since we have%
\begin{align*}
\left(  \left(  \left(  4\cdots4\right)  \cdot\left(  4\cdots0\right)
\right)  \cdot\left(  4\cdots4\right)  \right)  \cdot\left(  4\cdots4\right)
&  =\left(  \left(  4\cdot\ast\right)  \cdot4\right)  \cdot4=\left(
4\cdot4\right)  \cdot4=4\cdot4=4,\\
\left(  \left(  \left(  4\cdots4\right)  \cdot\left(  0\cdots0\right)
\right)  \cdot\left(  4\cdots4\right)  \right)  \cdot\left(  4\cdots4\right)
&  =\left(  \left(  4\cdot0\right)  \cdot4\right)  \cdot4=\left(
1\cdot4\right)  \cdot4=2\cdot4=3,
\end{align*}
and therefore in case $l_{2}>1$ we have%
\[
l_{3}=\max\Big\{i~\Big|~t^{\mathbf{A}}%
\bigl(%
\overset{l_{2}}{\overbrace{\underset{i}{\underbrace{4,\ldots,4}},0,\ldots0}%
},4,\ldots,4%
\bigr)%
=3\Big\}.
\]
Now we prove the inclusion $\sigma_{n} \left(  \mathbf{A}\right)
\supseteq\Sigma_{n}$, i.e. the fact that the numbers $l_{1},l_{2},l_{3}$
determine the term function $t^{\mathbf{A}}$. First we observe that
$\mathbf{A}$ satisfies the identity $x\left(  yz\right)  \approx xy$, and from
this we conclude by induction that%
\[
\mathbf{A} \models \eq{\en{1} \left(  t\right)}{ \en{1} \left( \operatorname{left}
\left(  t\right) \right) \cdot x_{l_{1}+1}}.
\]
Applying this identity to the left factor of $t$ we obtain%
\[
\mathbf{A}\models \eq{\en{1} \left( t \right)}{\left(
\en{1} \left( \operatorname{left}^{2}\left(  t\right)\right)  \cdot x_{l_{2}+1}\right)  \cdot
x_{l_{1}+1}}.
\]
Let us repeat this procedure until the left factor becomes the single variable
$x_{1}$. Suppose this happens after $s$ steps, i.e. $1=l_{s}<l_{s-1}%
<\ldots<l_{2}<l_{1}$. Then we have%
\[
\mathbf{A}\models\varepsilon_{1}\left(  t\right)  \approx\left(  \left(
\cdots\left(  \left(  x_{1}\cdot x_{l_{s}+1}\right)  \cdot x_{l_{s-1}%
+1}\right)  \cdots\right)  \cdot x_{l_{2}+1}\right)  \cdot x_{l_{1}+1.}%
\]
This already shows that $t^{\mathbf{A}}$ is determined by the numbers
$l_{1},l_{2},\ldots,l_{s}$. We have to show that actually the first three of
these numbers are sufficient. If $s\leq3$ then we have nothing to prove, and
if $s\geq4$, then using (\ref{nulla}) and the multiplication table we get the
following formula for $t^{\mathbf{A}}$:\medskip\\%
\parbox{\textwidth}{\[
\begin{array}{rcl}
t^{\mathbf{A}} \left(  x_{1},\ldots,x_{n}\right) &
= & \left(  \left(  \cdots\left(  \left(
x_{1}\cdot x_{l_{s}+1}\right)  \cdot x_{l_{s-1}+1}\right)  \cdots\right)
\cdot x_{l_{2}+1}\right)  \cdot x_{l_{1}+1} \smallskip \\
& = & \left\{  \!\!\!%
\begin{tabular}
[c]{ll}%
$0,$ & if $x_{1}=0;$\\
$1,$ & if $x_{1}\neq0=x_{l_{1}+1};$\\
$2,$ & if $x_{1},x_{l_{1}+1}\neq0=x_{l_{2}+1};$\\
$3,$ & if $x_{1},x_{l_{1}+1},x_{l_{2}+1}\neq0=x_{l_{3}+1};$\\
$4,$ & if $x_{1},x_{l_{1}+1},x_{l_{2}+1},x_{l_{3}+1}\neq0$.
\end{tabular}
\ \ \right.
\end{array}
\]
\eproof}
The next proposition shows that every finally associative spectrum appears as the fine spectrum of a finite groupoid.

\bpro\label{PropFinallyAssoc}
Let ${\bf A} \in \Alg{\omega}$ be a $p$-ary groupoid with
\[
 \exists n \in \N, n \ge 2 : \ \aS{{\bf A}}{n} = 1.
\]
Then there exists a finite groupoid ${\bf B} \in \Alg{\omega}$ with $\fS{{\bf B}} = \fS{\bf A}$.
\epro
\bproof
Let ${\bf A} \in \Alg{\omega}$ be a groupoid with the above property and denote by
$\Sigma := \fS{\bf A}$ the fine spectrum of ${\bf A}$. We know from Theorem~\ref{TheoremGenAssLaw} that
\[
 \forall m \in \N, m \ge n : \ \Sigma_m = \B{m}{p} \times \B{m}{p}
\]
holds.
And by Theorem~\ref{TheoremFS} we know that $\fS{\Tp / \Sigma} = \Sigma$ holds.
Define ${\bf B} = \pair{B}{\omega^{\bf B}}$ as
\[
 B := \rset{\Ec{t}_{\Sigma}}{t \in \B{k}{p}, k < n} \cup \set{*}
\]
with the operation
\[
 \function{\omega^{\bf B}}{B^p}{B}{\left(\Ec{t_1}_{\Sigma}, \ldots, \Ec{t_p}_{\Sigma}\right)}
 {\begin{cases}
   \Ec{\omega t_1 \ldots t_p}_{\Sigma} & \textrm{if } \occN{\omega t_1 \ldots t_p} < n \\
   * & \textrm{otherwise}
  \end{cases}}
\]
and $\omega^{\bf B} \left( b_1, \ldots, b_p \right) := *$ if one of the arguments is $*$.

We have to show that $\fS{{\bf B}} = \fS{\Tp / \Sigma}$ holds. This is pretty clear because
${\bf B}$ is nearly the same as $\Tp / \Sigma$. The only difference is that the equivalence classes containing all
bracketings of one size $m \ge n$ are equalized to $*$.
\eproof

\section{Open problems}

In conclusion, we formulate a few problems:

\begin{enumerate}[1.]
 \item Another idea to unterstand the lattice $\FSB$ is to translate constructions for groupoids into constructions in $\FSB$ and vice versa.
       A very simple example of this is the direct product $\prod$ and the meat $\bigwedge$:

       Let ${\bf A_i} \in \Alg{\omega}$ for $i \in I$ (arbitrary index set). Then we have:
       \[
        \fS{\prod_{i \in I} {\bf A_i}} = \bigwedge_{i \in I} \fS{{\bf A_i}}.
       \]
       Are there other correspondences between certain constructions, e.g. the join $\bigvee$ in $\FSB$?
 \item We have studied the {\sc Galois}-closed sets $\IdBrack \K$ for any $\K \subseteq \Alg{\omega}$.
       What is the analogon of a variety, i.e. what are the {\sc Galois}-closed sets $\ModBrack \Sigma$ on the groupoid side?
 \item What additional properties have fine spectra of finite algebras? Prove or disprove that
       the following condition is sufficient in order to realize a fine spectrum $\sigma$ on a finite groupoid:
       \[
        \exists N\in\N \ \forall n \in \N, n > N : \ \sigma_{n}=\delta_{n-1}\left( \sigma_{n-1}\right).
       \]

\end{enumerate}

\end{document}